\documentclass{amsart}
\usepackage{latexsym, amssymb, amsmath}

\def\RP{\mathbb {RP}}

\newtheorem{thm}{Theorem}[section]

\newtheorem{cor}[thm]{Corollary}
\newtheorem{prop}[thm]{Proposition}
\theoremstyle{remark}
\newtheorem{rmk}[thm]{Remark}
\theoremstyle{definition}
\newtheorem{defi}[thm]{Definition}
\newtheorem*{thmA}{Theorem~A}
\newtheorem*{thmB}{Theorem~B}

\title{Computations of the orbifold Yamabe invariant}

\author{Kazuo Akutagawa${}^*$} 
\email{akutagawa@math.is.tohoku.ac.jp} 
\address{Division of Mathematics, Graduate School of Information Sciences, Tohoku University, 
Sendai 980-8579, Japan} 

\thanks{${}^*$\ 
supported in part by the Grants-in-Aid for Scientific Research (C), 
Japan Society for the Promotion of Science, No.~21540059.} 

\date{September, 2010.} 
\begin{document} 
\maketitle 
\markboth{Computations of the orbifold Yamabe invariant} 
{Kazuo Akutagawa} 

\begin{abstract}

We consider the Yamabe invariant of a compact orbifold with finitely many singular points. 
We prove a fundamental inequality for the estimate of the invariant from above, 
which also includes a criterion for the non-positivity of it.  
Moreover, we give a sufficient condition for the equality in the inequality. 
In order to prove it, we also solve the orbifold Yamabe problem under a certain condition. 
We use these results to give some exact computations 
of the Yamabe invariant of compact orbifolds. 
\end{abstract}

\section{\bf Introduction and Main Results}\label{intro} 

There is a natural differential-topological invariant, 
called the {\it Yamabe invariant}, 
which arises from a variational problem for the functional $E$ below 
on a given compact smooth $n$-manifold $M$ (without boundary) of dimension $n \geq 3$. 
It is well known that a Riemannian metric on $M$ 
is {\it Einstein} if and only if 
it is a critical point of the normalized Einstein-Hilbert functional $E$ 
on the space $\mathcal{M}(M)$ of all Riemannian metrics on $M$ 
$$E : \mathcal{M}(M) \rightarrow \mathbb{R},\quad 
g \mapsto E(g) := \frac{\int_MR_g d\mu_g}{\textrm{Vol}_g(M)^{(n-2)/n}}. $$ 
Here, $R_g, d\mu_g$ and $\textrm{Vol}_g(M)$ denote respectively 
the scalar curvature, the volume element of $g$ and the volume of $(M, g)$. 
Because the restriction of $E$ to any conformal class 
$$C = [g] := \{ \textrm{e}^{2f}\cdot g\ |\ f \in C^{\infty}(M) \}$$ 
is bounded from below, 
we can consider the following conformal invariant 
(called the {\it Yamabe constant} of $(M, C)$) 
$$ 
Y(M, C) := \inf_{\tilde{g} \in C} E(\tilde{g}).  
$$ 
A remarkable theorem \cite{Yamabe, Trudinger, Aubin, Schoen-1, ScYa-3} 
(cf.~\cite{Aubin-Book, BrLi, LePa, Schoen-3, ScYa-Book}) of 
Yamabe, Trudinger, Aubin, and Schoen asserts that 
each conformal class $C$ contains a minimizer $\hat{g}$ of $E|_C$, 
called a {\it Yamabe metric} (or a {\it solution of the Yamabe problem}), 
which is of constant scalar curvature 
$$ 
R_{\hat{g}} = Y(M, C)\cdot \textrm{Vol}_{\hat{g}}(M)^{-2/n}. 
$$ 
The study of the second variation of $E$ 
done in \cite{Koiso, Schoen-3} (cf.~\cite{Besse-Book}) 
leads naturally to the definition of 
the following differential-topological invariant 
$$ 
Y(M) := \sup_{C \in \mathcal{C}(M)} \inf_{g \in C} E(g) = \sup_{C \in \mathcal{C}(M)} Y(M, C),  
$$ 
where $\mathcal{C}(M)$ denotes the space of all conformal classes on $M$. 
This invariant is called the {\it Yamabe invariant} (or {\it $\sigma$-invariant}) of $M$ 
and it was introduced independently by 
O.~Kobayashi~\cite{Kobayashi-1} and Schoen~\cite{Schoen-2} 
(see also \cite{Kobayashi-2, Schoen-3}). 

In the study of Yamabe invariant, 
with certain geometric non-collapsing assumptions, 
we will often encounter {\it Riemannian orbifolds} (or {\it Riemannian multi-folds} more generally) 
as the limit spaces for sequences of Yamabe metrics (cf.~\cite{Ak, TiVi, Vi}).  
For a compact $n$-orbifold $M$ with an orbifold metric $g$, 
one can also define the corresponding Yamabe constant $Y(M, [g]_{orb})$ and Yamabe invariant $Y^{orb}(M)$ 
(see Section~2 or \cite{AkBo-2} for details). 
Let $M_1$ and $M_2$ be compact $n$-orbifolds with same number of finite singularities 
$\{ \check{p}_1, \cdots, \check{p}_{\ell} \}$ 
and $\{ \check{q}_1, \cdots, \check{q}_{\ell} \}$ 
respectively. 
Assume that each corresponding singularities $\check{p}_j$ and $\check{q}_j$ 
have a same structure group $\Gamma_j ( < O(n) )$. 
For each $j$, 
let $B(\check{p}_j) ( \subset M_1 )$ and $B(\check{q}_j) ( \subset M_2 )$ denote respectively 
open geodesic balls of sufficiently small radiuses centered at $\check{p}_j$ and $\check{q}_j$ 
with fixed reference orbifold metrics. 
Then, the boundaries of these two balls can be naturally identified by a canonical diffeomorphism. 
Let 
$$ 
N := \big{(} M_1 - \sqcup_{j = 1}^{\ell} B(\check{p}_j) \big{)} \cup_Z 
\big{(} M_2 - \sqcup_{j = 1}^{\ell} B(\check{q}_j) \big{)}
$$ 
be the sum of $M_1 - \sqcup_{j = 1}^{\ell} B(\check{p}_j)$ and 
$M_2 - \sqcup_{j = 1}^{\ell} B(\check{q}_j)$ 
along their common boundary 
$Z := \partial \big{(} \sqcup_{j = 1}^{\ell} B(\check{p}_j) \big{)} 
= \partial \big{(} \sqcup_{j = 1}^{\ell} B(\check{q}_j) \big{)}$. 
Note that $N$ has a canonical smooth structure as manifold. 
For simplicity, in Section~4, we will abbreviate the above decomposition as 
the generalized connected sum 
$$ 
N = M_1 \#_{\sqcup_{j=1}^{\ell} (S^{n-1}/\Gamma_j)} M_2. 
$$ 

One of main purposes of this paper is to prove the following fundamental inequality 
for the estimate of the orbifold Yamabe invariant from above 
and a sufficient condition for the equality in this inequality. 
The inequality also includes a criterion for the non-positivity of the invariant:    

\begin{thmA} 
Under the above understandings, 
assume that 
$$ 
Y(N) \leq 0~({\rm resp.}~< 0)\quad {\rm and }\quad Y^{orb}(M_2) > 0~({\rm resp.}~\geq 0).  
$$ 
Then, 
$$ 
Y^{orb}(M_1) \leq Y(N) \leq 0. 
$$ 
Moreover, 
if $M_1$ admits an orbifold metric $\check{g}$ 
of constant scalar curvature satisfying $E(\check{g}) = Y(N)$, 
then 
$$ 
Y^{orb}(M_1) = Y(M_1, [\check{g}]_{orb}) = Y(N) \leq 0. 
$$ 
\end{thmA} 

On the computation of Yamabe invariants for {\it smooth} manifolds, 
a first remarkable result is the following proved by Aubin~\cite{Aubin} (cf.~\cite{Aubin-Book}) :  
$$ 
Y(M, C) \leq Y(S^n, [g_0]) = E(g_0)~\Big{(} = n(n-1) {\rm Vol}_{g_0}(S^n)^{2/n}~\Big{)}
$$ 
for any $C \in \mathcal{C}(M)$, 
where $g_0$ is the standard metric of constant curvature one on the standard $n$-sphere $S^n$. 
This implies both the universal estimate for $Y(M)$ from above and the computation of $Y(S^n)$  
$$ 
Y(M) \leq Y(S^n) = n(n-1) {\rm Vol}_{g_0}(S^n)^{2/n}. 
$$ 
Kobayashi~\cite{Kobayashi-1, Kobayashi-2} and Schoen~\cite{Schoen-3} proved that 
$$ 
Y(S^{n-1} \times S^1) = Y(S^n).  
$$ 
Kobayashi also gave two kind of proof for it (see \cite{AkFlPe} for the third one), 
one \cite{Kobayashi-2} of them especially is based on the following important inequality, 
called {\it Kobayashi's inequality}:  
$$ 
Y(M_1^n \# M_2^n) \geq 
\begin{cases}\ 
- \bigl(|Y(M_1^n)|^{n/2}+|Y(M_2^n)|^{n/2}\bigr)^{2/n}\ \cdots\ 
Y(M_1^n), Y(M_2^n) \leq 0, \\
\qquad \min \{ Y(M_1^n), Y(M_2^n) \}\ \cdots\ \text{otherwise} 
\end{cases} 
$$ 
for any two compact $n$-manifolds $M_1, M_2$. 
This has been extended to some useful surgery theorems \cite{AmDaHu, Petean-1, PeYu}. 
On the other hand, 
some classification theorems for manifolds with positive scalar curvature metric  
\cite{GrLa-1, GrLa-2, ScYa-1, ScYa-2, Stolz} lead to many examples of manifolds 
with zero (or non-positive) Yamabe invariant, for instance, $Y(T^n) = 0$ for the $n$-torus $T^n$ 
(see \cite{Petean-2} for further development). 

In 1995, LeBrun~\cite{Le-1} computed the Yamabe invariants of smooth compact quotients of complex-hyperbolic $2$-space,  
which was the first example of manifolds with negative Yamabe invariant. 
He and collaborators~\cite{GuLe, IsLe, Le-2, Le-3, Le-4} also computed the Yamabe invariants 
for a large class of $4$-manifolds, including K\"ahler surfaces 
$X$ with either $Y(X) < 0$ or $0 < Y(X) < Y(S^4)$ 
(see~\cite{AkNe, Anderson, BrNe} for $3$-manifolds $M^3$ with either $Y(M^3) < 0$ or $0 < Y(M^3) < Y(S^3)$~).  
In particular, for any minimal complex surface of general type $X$, 
he~\cite{Le-2} computed its Yamabe invariant $Y(X)$ to be   
$$ 
Y(X) = - 4 \sqrt{2} \pi \sqrt{2 \chi(X) + 3 \tau(X)} < 0, 
$$ 
where $\chi(X)$ and $\tau(X)$ are respectively the Euler characteristic and signature of $X$.  
Moreover, if $X$ contains $(- 2)$-curves, there exist a sequence of metrics $\{ g_i \}_i$ on $X$  
and a K\"ahler-Einstein orbifold metric $\check{g}$ on the canonical model $X_{can}$ of $X$ 
such that 
$$ 
\lim_{i \to \infty} Y(X, [g_i]) = Y(X),\qquad \lim_{i \to \infty} d_{GH} \big{(} (X, g_i), (X_{can}, \check{g}) \big{)} = 0. 
$$ 
Here, $d_{GH}$ denotes the Gromov-Hausdorff distance.  
This result suggests naturally the following question : 
``~Can one describe rigorously the above fact 
in terms of $Y(X_{can}, [\check{g}]_{orb})$ and $Y^{orb}(X_{can})$~?~'' 

The other of main purposes of this paper is to answer it. 

\begin{thmB} 
Under the above settings, 
the following holds 
$$  
Y^{orb}(X_{can}) = Y(X_{can}, [\check{g}]_{orb}) = Y(X).    
$$ 
\end{thmB} 

In Section~2, we recall the definition on the orbifold Yamabe invariant from \cite{AkBo-2}  
and explain briefly some terminologies. 
For the proof of Theorem~A, 
we also recall some necessary terminologies and basic results 
on the Yamabe invariant of cylindrical manifolds~\cite{AkBo-1}.  
Applying these results to the orbifold Yamabe invariant, 
we prove the first assertion of Theorem~A. 
In Section~3, for the proof of the second assertion in Theorem~A, 
we consider the existence problem of minimizers for the functional $E$ 
on compact conformal orbifolds, that is, the {\it orbifold Yamabe problem}. 
Under a certain condition, 
we solve this problem. 
Using the solution, we can prove the second assertion. 
In Section~4, we give two more typical exact computations of the orbifold Yamabe invariant 
besides the proof of Theorem~B.   

\noindent 
{\bf Acknowledgements.} 
The author would like to express his sincere gratitude 
to Nobuhiro Honda and Jeff Viaclovsky for helpful discussions 
on singularities of complex surfaces 
and on the orbifold Yamabe invariant respectively. 
He also would like to thank Claude LeBrun for useful comments.

\section{\bf The orbifold Yamabe invariant}\label{OYI} 

For the sake of self-containedness, 
we first recall the definition of orbifolds with finitely many singular points 
which we discuss here \cite{AkBo-2}. 

\begin{defi} 
Let $M$ be a locally compact Hausdorff space. 
We say that $M$ is an $n$-{\it orbifold with singularities} 
$$ 
\Sigma_{\Gamma} = \{(\check{p}_1, \Gamma_1), \cdots, (\check{p}_{\ell}, \Gamma_{\ell})\} 
$$ 
if the following conditions are satisfied: \\ 
\quad (1) $\Sigma := \{\check{p}_1, \cdots, \check{p}_{\ell}\} \subset M$, 
and $M - \Sigma$ is a smooth $n$-manifold. \\ 
\quad (2) $\Gamma := \{\Gamma_1, \cdots, \Gamma_{\ell}\}$ is a collection of non-trivial finite subgroups $\Gamma_j$ of $O(n)$, 
each of which acts freely on $\mathbb{R}^n - \{{\bf 0}\}$.  \\ 
\quad (3) For each $\check{p}_j$, there exist its open neighborhood $U_j$ 
and a homeomorphism $\varphi_j : U_j \rightarrow \mathbb{B}_{\tau_j}({\bf 0})/\Gamma_j$ for some $\tau_j > 0$ such that 
$$ 
\varphi_j : U_j - \{\check{p}_j\} \longrightarrow \big{(}\mathbb{B}_{\tau_j}({\bf 0}) - \{{\bf 0}\}\big{)}/\Gamma_j 
$$ 
is a diffeomorphism. 
Here, $\mathbb{B}_{\tau_j}({\bf 0}) := \{ x = (x^1, \cdots, x^n) \in \mathbb{R}^n~|~|x| < \tau_j \}$. 
\end{defi} 

We refer to the pair $(\check{p}_j, \Gamma_j)$ as a {\it singular point with the structure group} $\Gamma_j$ 
and the pair $(U_j, \varphi_j)$ as a {\it local uniformization}. 
To simplify the presentation, we assume, without particular mention, 
that an orbifold $M$ has only one singularity, i.e., $\Sigma_{\Gamma} = \{(\check{p}, \Gamma)\}$. 
Let $\varphi : U \rightarrow \mathbb{B}_{\tau}({\bf 0})/\Gamma$ be a local uniformization 
and $\pi : \mathbb{B}_{\tau}({\bf 0}) \rightarrow \mathbb{B}_{\tau}({\bf 0})/\Gamma$ the canonical projection. 
We also always assume that $M$ is compact. 

\begin{defi} 
(1)\ \ A Riemannian metric $g \in \mathcal{M}(M - \{p\})$ is an {\it orbifold metric} 
if there exists a $\Gamma$-invariant smooth metric $\hat{g}$ on the ball $\mathbb{B}_{\tau}({\bf 0})$ 
such that 
$(\varphi^{-1} \circ \pi)^{\ast} g = \hat{g}$ on $\mathbb{B}_{\tau}({\bf 0}) - \{{\bf 0}\}$. 
We denote by $\mathcal{M}^{orb}(M)$ the space of all orbifold metrics on $M$. 
In the case when 
$\Sigma_{\Gamma} = \{(\check{p}_1, \Gamma_1), \cdots, (\check{p}_{\ell}, \Gamma_{\ell})\}$, 
the space of all orbifold metrics is defined similarly. \\ 
\quad (2)\ \ For an orbifold metric $g \in \mathcal{M}^{orb}(M)$, 
its {\it orbifold conformal class} $[g]_{orb}$ is defined by 
\begin{align*} 
[g]_{orb} &:= [g] \cap \mathcal{M}^{orb}(M) \\ 
&\ = \{ e^{2f}\cdot g~|~f \in C^0(M) \cap C^{\infty}(M - \{\check{p}\}), 
(\varphi^{-1}\circ \pi)^{\ast}f \in C^{\infty}(\mathbb{B}_{\tau}({\bf 0})) \}.
\end{align*} 
We denote by $\mathcal{C}^{orb}(M)$ the space of all orbifold conformal classes. 
\end{defi} 

As in the smooth case, consider the normalized Einstein-Hilbert functional 
$$ 
E : \mathcal{M}^{orb}(M) \rightarrow \mathbb{R},\quad 
g \mapsto \frac{\int_MR_g d\mu_g}{\textrm{Vol}_g(M)^{(n-2)/n}}. 
$$  
Since the singularity has codimension at least three, 
Stokes' theorem and Gauss' divergence theorem still hold over Riemannian orbifolds. 
Hence, $\check{g}$ is a critical point of $E$ on ${\mathcal{M}^{orb}(M)}$ 
if and only if $\check{g}$ is an Einstein orbifold metric. 
Then, one can define naturally the definition of the orbifold Yamabe invariant. 

\begin{defi} 
For a conformal orbifold $(M, [g]_{orb})$, 
its {\it Yamabe constant} $Y(M, [g]_{orb})$ is defined by 
$$ 
Y(M, [g]_{orb}) := \inf_{\tilde{g} \in [g]_{orb}} E(\tilde{g}). 
$$ 
Moreover, the {\it orbifold Yamabe invariant} $Y^{orb}(M)$ of $M$ is also defined by 
$$ 
Y^{orb}(M) := \sup_{[g]_{orb} \in \mathcal{C}^{orb}(M)}Y(M, [g]_{orb}). 
$$ 
\end{defi} 

Before we explain some necessary terminologies 
on the Yamabe invariant of cylindrical manifolds, 
we give two comments on orbifolds with {\it positive} orbifold Yamabe invariant. 

\begin{rmk} 
Let $(X, g)$ be a hyperK\"ahler asymptotically locally Euclidean (abbreviated to {\it ALE} ) $4$-manifold 
constructed in \cite{Kr} (cf.~\cite{Na}), 
where $X$ is the minimal resolution of the quotient space $\mathbb{C}^2/\Gamma$  
for a non-trivial finite subgroup $\Gamma$ of $SU(2)$. 
Then, $(X, g)$ has a smooth conformal compactification $(\hat{X} := X \sqcup \{p_{\infty}\}, \hat{g})$ 
with singularity $\{(p_{\infty}, \Gamma)\}$ \cite{ChLeWe, Vi}, 
which has a positive Yamabe constant $Y(\hat{X}, [\hat{g}]_{orb}) > 0$. 
In \cite[Theorem~1.3]{Vi}, Viaclovsky has proved the following: \\ 
\quad (1)\ \ The orbifold Yamabe problem on $(\hat{X}, [\hat{g}]_{orb})$ has no solution. 
This implies that the orbifold Yamabe problem is not always solvable (see Section~3 for the solvability), 
in contrast with the case for smooth compact conformal manifolds. \\ 
\quad (2)\ \ He computed the orbifold Yamabe invariant of $\hat{X}$ as 
$$ 
Y^{orb}(\hat{X}) = Y(\hat{X}, [\hat{g}]_{orb}) = Y(S^4)/|\Gamma|^{1/2}. 
$$ 
However, similarly to the case for smooth compact manifolds, 
there is not much exact computations of positive orbifold Yamabe invariants at present. 
  
In the proof of both (1) and (2), one of key points is the following estimate, 
called {\it refined Aubin's inequality} \cite[Theorem~B]{AkBo-2}
$$ 
Y(M, [g]_{orb}) \leq Y^{orb}(M) \leq \min_{1 \leq j \leq \ell} \frac{Y(S^n)}{|\Gamma_j|^{2/n}} 
$$ 
for any compact Riemannian $n$-orbifold $(M, g)$ with singularities 
$\{(\check{p}_1, \Gamma_1), \cdots, (\check{p}_{\ell}, \Gamma_{\ell})\}$. 
This inequality is also crucial to give a sufficient condition 
for the solvability of the orbifold Yamabe problem in Section~3.  
\end{rmk} 

\begin{defi} 
Let $X$ be an open $n$-manifold {\it with tame ends}, i.e., 
it is diffeomorphic to $W \cup_Z (Z \times [0, \infty))$, 
where $W (\subset X)$ is a relatively compact open submanifold with 
boundary $\partial W =: Z \cong Z \times \{0\}$ (possibly finitely many connected component). 
For a fixed $h \in \mathcal{M}(Z)$, 
a complete Riemannian metric $\bar{g}$ on $X$ is called a {\it cylindrical metric modeled by} $(Z, h)$ 
if there exists a global coordinate function $t$ on $Z \times [0, \infty)$ such that 
$\bar{g}|_{Z \times [1, \infty)}$ is the product metric $\bar{g}(z, t) = h(z) + dt^2~(~(z, t) \in Z \times [1, \infty)~)$ 
(see Figure~1). 
Each pair $(X, \bar{g})$ is called a {\it cylindrical manifold} and $h$ a {\it slice metric}. 
We denote by $\mathcal{M}^{h\textrm{-}cyl}(X)$ the space of all cylindrical metrics on $X$ modeled by $(Z, h)$. 
\end{defi} 
\quad \\ 
\unitlength 0.1in
\begin{picture}(42.02,11.00)(1.60,-11.60)
%
\special{pn 8}%
\special{ar 1044 610 884 550  0.9576225 5.3666843}%
%
\special{pn 8}%
\special{pa 1502 140}%
\special{pa 1532 152}%
\special{pa 1561 165}%
\special{pa 1591 177}%
\special{pa 1620 190}%
\special{pa 1649 204}%
\special{pa 1678 218}%
\special{pa 1706 232}%
\special{pa 1734 248}%
\special{pa 1762 264}%
\special{pa 1789 280}%
\special{pa 1817 296}%
\special{pa 1845 311}%
\special{pa 1874 325}%
\special{pa 1903 338}%
\special{pa 1933 350}%
\special{pa 1963 361}%
\special{pa 1994 372}%
\special{pa 2024 381}%
\special{pa 2055 390}%
\special{pa 2086 398}%
\special{pa 2118 405}%
\special{pa 2149 412}%
\special{pa 2180 418}%
\special{pa 2212 424}%
\special{pa 2243 428}%
\special{pa 2275 433}%
\special{pa 2307 436}%
\special{pa 2339 440}%
\special{pa 2370 443}%
\special{pa 2402 445}%
\special{pa 2434 447}%
\special{pa 2466 449}%
\special{pa 2498 450}%
\special{pa 2530 451}%
\special{pa 2563 452}%
\special{pa 2595 452}%
\special{pa 2627 452}%
\special{pa 2659 453}%
\special{pa 2691 452}%
\special{pa 2724 452}%
\special{pa 2756 452}%
\special{pa 2788 452}%
\special{pa 2820 451}%
\special{pa 2852 451}%
\special{pa 2885 450}%
\special{pa 2917 450}%
\special{pa 2949 450}%
\special{pa 2981 450}%
\special{pa 3013 449}%
\special{pa 3045 449}%
\special{pa 3077 449}%
\special{pa 3109 449}%
\special{pa 3141 449}%
\special{pa 3173 449}%
\special{pa 3205 450}%
\special{pa 3237 450}%
\special{pa 3269 450}%
\special{pa 3301 450}%
\special{pa 3333 450}%
\special{pa 3365 450}%
\special{pa 3397 450}%
\special{pa 3429 450}%
\special{pa 3460 450}%
\special{pa 3492 451}%
\special{pa 3524 451}%
\special{pa 3556 451}%
\special{pa 3588 451}%
\special{pa 3620 451}%
\special{pa 3652 451}%
\special{pa 3684 451}%
\special{pa 3716 451}%
\special{pa 3748 451}%
\special{pa 3780 451}%
\special{pa 3812 451}%
\special{pa 3844 451}%
\special{pa 3876 451}%
\special{pa 3908 451}%
\special{pa 3940 451}%
\special{pa 3972 451}%
\special{pa 4004 451}%
\special{pa 4036 451}%
\special{pa 4068 451}%
\special{pa 4100 450}%
\special{pa 4132 450}%
\special{pa 4165 450}%
\special{pa 4197 450}%
\special{pa 4229 450}%
\special{pa 4261 450}%
\special{pa 4293 450}%
\special{pa 4325 450}%
\special{pa 4357 450}%
\special{pa 4362 450}%
\special{sp}%
%
\special{pn 8}%
\special{pa 1439 1110}%
\special{pa 1467 1095}%
\special{pa 1495 1080}%
\special{pa 1524 1064}%
\special{pa 1552 1049}%
\special{pa 1580 1035}%
\special{pa 1609 1020}%
\special{pa 1637 1005}%
\special{pa 1666 991}%
\special{pa 1695 977}%
\special{pa 1724 963}%
\special{pa 1753 950}%
\special{pa 1782 937}%
\special{pa 1811 924}%
\special{pa 1841 912}%
\special{pa 1871 900}%
\special{pa 1900 889}%
\special{pa 1931 878}%
\special{pa 1961 867}%
\special{pa 1991 857}%
\special{pa 2021 847}%
\special{pa 2052 837}%
\special{pa 2083 828}%
\special{pa 2114 820}%
\special{pa 2145 812}%
\special{pa 2176 804}%
\special{pa 2207 797}%
\special{pa 2238 790}%
\special{pa 2270 784}%
\special{pa 2301 778}%
\special{pa 2333 772}%
\special{pa 2365 767}%
\special{pa 2396 762}%
\special{pa 2428 758}%
\special{pa 2460 754}%
\special{pa 2492 750}%
\special{pa 2523 746}%
\special{pa 2555 743}%
\special{pa 2587 741}%
\special{pa 2619 738}%
\special{pa 2651 736}%
\special{pa 2683 734}%
\special{pa 2715 732}%
\special{pa 2747 731}%
\special{pa 2779 730}%
\special{pa 2811 729}%
\special{pa 2843 728}%
\special{pa 2875 728}%
\special{pa 2907 727}%
\special{pa 2939 727}%
\special{pa 2971 727}%
\special{pa 3003 727}%
\special{pa 3035 727}%
\special{pa 3067 727}%
\special{pa 3099 727}%
\special{pa 3131 728}%
\special{pa 3163 728}%
\special{pa 3195 729}%
\special{pa 3227 729}%
\special{pa 3259 729}%
\special{pa 3291 730}%
\special{pa 3323 730}%
\special{pa 3355 731}%
\special{pa 3388 731}%
\special{pa 3420 731}%
\special{pa 3452 731}%
\special{pa 3484 732}%
\special{pa 3516 732}%
\special{pa 3548 732}%
\special{pa 3580 732}%
\special{pa 3612 732}%
\special{pa 3644 732}%
\special{pa 3676 732}%
\special{pa 3708 732}%
\special{pa 3740 732}%
\special{pa 3772 732}%
\special{pa 3804 732}%
\special{pa 3836 732}%
\special{pa 3868 732}%
\special{pa 3900 732}%
\special{pa 3932 732}%
\special{pa 3964 732}%
\special{pa 3996 732}%
\special{pa 4028 732}%
\special{pa 4060 732}%
\special{pa 4092 732}%
\special{pa 4124 731}%
\special{pa 4156 731}%
\special{pa 4188 731}%
\special{pa 4220 731}%
\special{pa 4251 731}%
\special{pa 4283 730}%
\special{pa 4315 730}%
\special{pa 4347 730}%
\special{pa 4362 730}%
\special{sp}%
%
\special{pn 8}%
\special{pa 1225 576}%
\special{pa 1241 604}%
\special{pa 1252 634}%
\special{pa 1255 666}%
\special{pa 1252 697}%
\special{pa 1242 728}%
\special{pa 1225 754}%
\special{pa 1203 778}%
\special{pa 1178 798}%
\special{pa 1151 814}%
\special{pa 1122 828}%
\special{pa 1091 839}%
\special{pa 1060 846}%
\special{pa 1029 852}%
\special{pa 997 855}%
\special{pa 965 856}%
\special{pa 933 856}%
\special{pa 901 854}%
\special{pa 870 848}%
\special{pa 838 842}%
\special{pa 807 833}%
\special{pa 777 823}%
\special{pa 748 811}%
\special{pa 719 796}%
\special{pa 693 778}%
\special{pa 667 759}%
\special{pa 644 736}%
\special{pa 624 712}%
\special{pa 608 684}%
\special{pa 598 654}%
\special{pa 594 622}%
\special{pa 596 590}%
\special{pa 607 560}%
\special{pa 624 533}%
\special{pa 646 510}%
\special{pa 651 505}%
\special{sp}%
%
\special{pn 8}%
\special{pa 586 640}%
\special{pa 610 619}%
\special{pa 638 603}%
\special{pa 667 590}%
\special{pa 697 579}%
\special{pa 728 572}%
\special{pa 760 567}%
\special{pa 791 562}%
\special{pa 823 560}%
\special{pa 855 560}%
\special{pa 887 560}%
\special{pa 919 562}%
\special{pa 951 566}%
\special{pa 983 570}%
\special{pa 1014 577}%
\special{pa 1045 585}%
\special{pa 1076 594}%
\special{pa 1106 605}%
\special{pa 1135 617}%
\special{pa 1164 632}%
\special{pa 1191 648}%
\special{pa 1217 668}%
\special{pa 1239 690}%
\special{pa 1256 711}%
\special{sp}%
%
\special{pn 8}%
\special{ar 1866 620 124 280  1.4056476 4.8723021}%
%
\special{pn 8}%
\special{ar 1876 630 146 280  4.7123890 4.7687270}%
\special{ar 1876 630 146 280  4.9377411 4.9940791}%
\special{ar 1876 630 146 280  5.1630932 5.2194312}%
\special{ar 1876 630 146 280  5.3884453 5.4447833}%
\special{ar 1876 630 146 280  5.6137974 5.6701355}%
\special{ar 1876 630 146 280  5.8391495 5.8954876}%
\special{ar 1876 630 146 280  6.0645017 6.1208397}%
\special{ar 1876 630 146 280  6.2898538 6.3461918}%
\special{ar 1876 630 146 280  6.5152059 6.5715439}%
\special{ar 1876 630 146 280  6.7405580 6.7968960}%
\special{ar 1876 630 146 280  6.9659101 7.0222481}%
\special{ar 1876 630 146 280  7.1912622 7.2476002}%
%
\special{pn 8}%
\special{ar 2490 610 93 150  4.9117263 5.0104917}%
\special{ar 2490 610 93 150  5.3067880 5.4055534}%
\special{ar 2490 610 93 150  5.7018497 5.8006152}%
\special{ar 2490 610 93 150  6.0969115 6.1956769}%
\special{ar 2490 610 93 150  6.4919732 6.5907386}%
\special{ar 2490 610 93 150  6.8870349 6.9858004}%
\special{ar 2490 610 93 150  7.2820967 7.3152542}%
%
\special{pn 8}%
\special{ar 2490 610 114 150  1.4767351 4.6161956}%
\put(28.3000,-3.7000){\makebox(0,0)[lb]{{\small $X - W \cong Z\times [0,\infty)$}}}%
%
\special{pn 8}%
\special{ar 3609 590 113 140  1.4767351 4.6153496}%
%
\special{pn 8}%
\special{ar 3590 600 93 150  4.9117263 5.0104917}%
\special{ar 3590 600 93 150  5.3067880 5.4055534}%
\special{ar 3590 600 93 150  5.7018497 5.8006152}%
\special{ar 3590 600 93 150  6.0969115 6.1956769}%
\special{ar 3590 600 93 150  6.4919732 6.5907386}%
\special{ar 3590 600 93 150  6.8870349 6.9858004}%
\special{ar 3590 600 93 150  7.2820967 7.3152542}%
\put(8.8000,-4.4000){\makebox(0,0)[lb]{$W$}}%
%
\special{pn 8}%
\special{pa 2490 920}%
\special{pa 4100 920}%
\special{fp}%
\special{sh 1}%
\special{pa 4100 920}%
\special{pa 4033 900}%
\special{pa 4047 920}%
\special{pa 4033 940}%
\special{pa 4100 920}%
\special{fp}%
\put(24.8000,-12.7000){\makebox(0,0)[lb]{{\small $\bar{g}(z,t) = h(z)+dt^2\ \text{on}\ Z\times [1,\infty)$}}}%
%
\special{pn 8}%
\special{pa 2490 920}%
\special{pa 2490 740}%
\special{fp}%
\end{picture}%
\\ 
\quad \\ 
\qquad \qquad \qquad Figure~1:~A cylindrical manifold $(X, \bar{g})$\\ 

For the definition of the Yamabe invariant on cylindrical manifolds, 
we first recall the following fact. 
On a compact Riemannian manifold $(M, g)$, 
the value of functional $E(\tilde{g})$ for conformal metric $\tilde{g} := u^{4/(n-2)}\cdot g \in [g]$ 
can be rewritten by 
$$ 
E(\widetilde{g}) = \frac{\int_M \big{(}\alpha_n|\nabla u|^2 + R_gu^2\big{)} d\mu_g}
{\Big{(}\int_M u^{2n/(n-2)} d\mu_g\Big{)}^{(n-2)/n}}\ \Big{(} =: Q_{(M, g)}(u)~\Big{)},  
\qquad \alpha_n := \frac{4(n-1)}{n-2} > 0.  
$$ 

\begin{defi} 
The {\it Yamabe constant} $Y(X, [\bar{g}])$ of a cylindrical manifold $(X, \bar{g})$ is defined by 
$$ 
Y(X, [\bar{g}]) := \inf_{u \in C^{\infty}_c(X), u \not\equiv 0} Q_{(X, \bar{g})}(u), 
$$ 
where $C^{\infty}_c(X)$ denotes the space of all smooth functions on $X$ with compact supports. 
Moreover, for a fixed $h \in \mathcal{M}(Z)$, 
the $h$-{\it cylindrical Yamabe invariant} $Y^{h\textrm{-}cyl}(X)$ of the open manifold $X$ with tame ends 
is also defined by 
$$ 
Y^{h\textrm{-}cyl}(X) := \sup_{\bar{g} \in \mathcal{M}^{h\textrm{-}cyl}(X)} Y(X, [\bar{g}]). 
$$ 
\end{defi} 

To simplify the presentation, we also assume, without particular mention, 
that each underlying manifold $X$ has only one connected tame end. 
In contrast with the case for compact manifolds, 
the constant $Y(X, [\bar{g}])$ is not always finite. 
For instance, if the scalar curvature $R_h$ of slice metric $h$ is negative on $Z$, 
then $Y(X, [\bar{g}]) = - \infty$. 
As a complete criterion for the finiteness of $Y(X, [\bar{g}])$, 
we have obtained the following \cite[Lemmas~2.7, 2.9]{AkBo-1}. 

\begin{prop} 
For $h \in \mathcal{M}(Z^{n-1})$, let $\mathcal{L}_h$ be the operator on $Z^{n-1}$ defined by 
$$ 
\mathcal{L}_h := - \frac{4(n-1)}{n-2}\Delta_h + R_h,  
$$ 
and $\lambda(\mathcal{L}_h)$ the first eigenvalue of $\mathcal{L}_h$. 
Then, we have the following on the Yamabe constant of a cylindrical manifold $(X, \bar{g})$ with slice metric $h$. \\ 
\quad $\bullet$\ \ If $\lambda(\mathcal{L}_h) < 0$, then $Y(X, [\bar{g}]) = - \infty$. \\ 
\quad $\bullet$\ \ If $\lambda(\mathcal{L}_h) \geq 0$, then $Y(X, [\bar{g}]) > - \infty$. \\ 
\quad $\bullet$\ \ If $\lambda(\mathcal{L}_h) = 0$, then $0 \geq Y(X, [\bar{g}]) > - \infty$. 
\end{prop}  

We also note that the notion of the $h$-cylindrical Yamabe invariant is an natural extension of the one 
of the orbifold Yamabe invariant \cite[Theorem~2.9]{AkBo-2}. 

\begin{prop} 
Let $M$ be a compact $n$-orbifold with singularity $\{(\check{p}, \Gamma)\}$ $($see {\rm Figure~2}$)$, 
and $h_0 \in \mathcal{M}(S^{n-1}/\Gamma)$ the standard metric of constant curvature one. 
Note that the open manifold $M - \{\check{p}\}$ is of one tame end 
and $\mathcal{M}^{h_0\textrm{-}cyl}(M - \{\check{p}\}) \not= \emptyset$. 
Then, 
$$ 
Y^{orb}(M) = Y^{h_0\textrm{-}cyl}(M - \{\check{p}\}). 
$$ 
\end{prop} 
\quad \\ 
\qquad \qquad \qquad \qquad \qquad \qquad \qquad 
\unitlength 0.1in
\begin{picture}(18.70,11.00)(34.90,-11.70)
%
\special{pn 8}%
\special{ar 5115 575 45 115  4.9177844 5.0677844}%
\special{ar 5115 575 45 115  5.5177844 5.6677844}%
\special{ar 5115 575 45 115  6.1177844 6.2677844}%
\special{ar 5115 575 45 115  6.7177844 6.8677844}%
%
\special{pn 8}%
\special{ar 5115 575 55 115  1.4730694 4.6217291}%
%
\special{pn 8}%
\special{ar 4255 620 765 550  0.9579761 5.3670332}%
%
\special{pn 8}%
\special{pa 4412 586}%
\special{pa 4425 615}%
\special{pa 4435 646}%
\special{pa 4438 677}%
\special{pa 4434 709}%
\special{pa 4425 740}%
\special{pa 4409 768}%
\special{pa 4389 792}%
\special{pa 4364 813}%
\special{pa 4338 830}%
\special{pa 4309 844}%
\special{pa 4278 854}%
\special{pa 4247 861}%
\special{pa 4215 865}%
\special{pa 4183 866}%
\special{pa 4151 866}%
\special{pa 4120 861}%
\special{pa 4088 855}%
\special{pa 4058 846}%
\special{pa 4028 834}%
\special{pa 3999 821}%
\special{pa 3972 804}%
\special{pa 3947 784}%
\special{pa 3923 763}%
\special{pa 3903 738}%
\special{pa 3886 711}%
\special{pa 3873 681}%
\special{pa 3867 650}%
\special{pa 3865 618}%
\special{pa 3870 587}%
\special{pa 3883 557}%
\special{pa 3901 531}%
\special{pa 3915 515}%
\special{sp}%
%
\special{pn 8}%
\special{pa 3859 650}%
\special{pa 3882 628}%
\special{pa 3909 610}%
\special{pa 3938 597}%
\special{pa 3967 585}%
\special{pa 3998 578}%
\special{pa 4030 573}%
\special{pa 4062 570}%
\special{pa 4094 569}%
\special{pa 4126 570}%
\special{pa 4158 573}%
\special{pa 4189 578}%
\special{pa 4221 585}%
\special{pa 4252 593}%
\special{pa 4282 604}%
\special{pa 4311 617}%
\special{pa 4340 630}%
\special{pa 4367 648}%
\special{pa 4392 667}%
\special{pa 4416 689}%
\special{pa 4434 715}%
\special{pa 4439 721}%
\special{sp}%
%
\special{pn 8}%
\special{pa 4678 160}%
\special{pa 5272 540}%
\special{fp}%
%
\special{pn 8}%
\special{pa 4651 1100}%
\special{pa 5263 570}%
\special{fp}%
%
\special{pn 8}%
\special{ar 4894 600 108 280  1.3994172 4.8775377}%
%
\special{pn 8}%
\special{pa 4894 330}%
\special{pa 4921 344}%
\special{pa 4942 368}%
\special{pa 4955 398}%
\special{pa 4966 427}%
\special{pa 4972 459}%
\special{pa 4978 490}%
\special{pa 4981 522}%
\special{pa 4984 554}%
\special{pa 4984 586}%
\special{pa 4984 618}%
\special{pa 4982 650}%
\special{pa 4978 682}%
\special{pa 4975 713}%
\special{pa 4969 745}%
\special{pa 4961 776}%
\special{pa 4952 806}%
\special{pa 4941 836}%
\special{pa 4928 865}%
\special{pa 4919 880}%
\special{sp -0.045}%
%
\special{pn 20}%
\special{sh 1}%
\special{ar 5270 560 10 10 0  6.28318530717959E+0000}%
\special{sh 1}%
\special{ar 5270 550 10 10 0  6.28318530717959E+0000}%
\put(53.6000,-5.7000){\makebox(0,0)[lb]{$\check{p}$}}%
\put(48.9000,-9.9000){\makebox(0,0)[lt]{$S^{n-1}/\Gamma$}}%
\end{picture}%
\\ 
\quad \\ 
\qquad \qquad \qquad \qquad \qquad \qquad \qquad \qquad \qquad Figure~2. \\ 

Now, we can state the key inequality for $h$-cylindrical Yamabe invariants, 
called {\it refined Kobayashi's inequality} \cite[Theorem~3.7]{AkBo-1}. 

\begin{thm} 
Let $N$ be a compact $n$-manifold and $Z$ a compact $(n-1)$-submanifold 
with trivial normal bundle. 
Assume that $M - Z$ has two connected components $W_1, W_2$. 
Let $X_1 := \overline{W}_1 \cup_Z (Z \times [0, \infty)), 
X_2 := \overline{W}_2 \cup_Z (Z \times [0, \infty))$ 
be the corresponding open $n$-manifolds with tame end $Z \times [0, \infty)$ 
$($see {\rm Figure~3}$)$. 
For any $h \in \mathcal{M}(Z)$, we have 
$$ 
Y(N) \geq 
{\small \begin{cases}\ 
- \big{(}Y^{h\text{-}cyl}(X_1)^{n/2}+|Y^{h\text{-}cyl}(X_2)|^{n/2}\big{)}^{2/n}\ \cdots\ 
if\ \ Y^{h\text{-}cyl}(X_1), Y^{h\text{-}cyl}(X_2) \leq 0, \\
\ \min \{ Y^{h\text{-}cyl}(X_1), Y^{h\text{-}cyl}(X_2) \}\ \cdots\ \text{otherwise}. \end{cases} }
$$ 
\end{thm}

\qquad \qquad \qquad $N = \overline{W}_1 \cup_Z \overline{W}_2$ \\ 
\quad \\ 
\input{fig-n5.tex} \\
\quad \\ 
\qquad \qquad \qquad \qquad \qquad \qquad \qquad \qquad \qquad Figure~3. \\ 

Theorem~2.9 implies immediately the following. 

\begin{cor} 
Under the same setting as in Theorem~2.9, 
assume that $Y(N) \leq 0\\ ({\rm resp.}~< 0)$ and 
$Y^{h\textrm{-}cyl}(X_2) > 0~({\rm resp.}~\geq 0)$.  
$($From Proposition~2.7, the positivity $Y^{h\textrm{-}cyl}(X_2) > 0$ 
implies automatically $\lambda(\mathcal{L}_h) > 0.)$ 
Then, we have 
$$ 
Y^{h\textrm{-}cyl}(X_1) \leq Y(N) \leq 0.  
$$ 
\end{cor} 

We can now prove the first assertion in Theorem~A. \\ 
\quad \\ 
{\it Proof of the first assertion in Theorem~A}.\ \ 
In Corollary~2.10, 
set $W_1 = M_1 - \sqcup_{j = 1}^{\ell} B(\check{p}_j)$,  
$W_2 = M_2 - \sqcup_{j = 1}^{\ell} B(\check{q}_j)$ 
and $h = h_0$ 
on $Z = \partial W_1 = \partial W_2 \cong \sqcup _{j=1}^{\ell}\big{(}S^{n-1}/\Gamma_j\big{)}$. 
Note that $X_1 = M_1 - \{\check{p}_1, \cdots, \check{p}_{\ell}\}$ and 
$X_2 = M_2 - \{\check{q}_1, \cdots, \check{q}_{\ell}\}$. 
Then, the first assertion follows directly from Proposition~2.8 and Corollary~2.10, that is, 
$$ 
\qquad \qquad \qquad \qquad \qquad Y^{orb}(M_1) = Y^{h_0\textrm{-}cyl}(X_1) \leq Y(N) \leq 0. 
\qquad \qquad \qquad \qquad \quad \ \ \square
$$ 

\begin{rmk} 
For given compact manifolds $N_1$ and $N_2$,  
we generally use Kobayashi's inequality in the case for 
computing (or estimating) $Y(N_1 \# N_2)$ 
by using the values of both $Y(N_1)$ and $Y(N_2)$. 
In contrast with this, 
the generalized connected sum of compact orbifolds 
is often ``prime'' as smooth manifold. 
Hence, the opposite usage of (refined) Kobayashi's inequality 
is also useful as Theorem~A.  
\end{rmk}

\section{\bf The orbifold Yamabe problem}\label{EC} 

In this section, we first prove the orbifold Yamabe problem 
under a certain condition. 

\begin{thm} 
Let $(M, g)$ be a compact Riemannian $n$-orbifold with singularities 
$\{(\check{p}_1, \Gamma_1), \cdots, (\check{p}_{\ell}, \Gamma_{\ell})\}$. 
Assume the following strict inequality$:$   
\begin{equation} 
Y(M, [g]_{orb}) <~\min_{1 \leq j \leq \ell} \frac{Y(S^n)}{|\Gamma_j|^{2/n}}.
\end{equation} 
Then, there exists a minimizer $\tilde{g} \in [g]_{orb}$ 
of the functional $E|_{[g]_{orb}}$ $($called an {\rm orbifold Yamabe metric}$)$  
such that the orbifold metric $\tilde{g}$ is of constant scalar curvature 
$R_{\tilde{g}} = Y(M, [g]_{orb})\cdot {\rm Vol}_{\tilde{g}}(M)^{-2/n}$. 
\end{thm} 
\noindent  
{\it Proof.}\ \ We use here the same notations as those in Definition~2.1. 
Without loss of generality, 
we may assume that $M$ has only one singularity $\{(\check{p}, \Gamma)\}$. 
The method adopting here for constructing approximate solutions 
is similar to the one in \cite[Theorem~5.2]{AkBo-1}. 
But, as background metric for getting both the uniform $C^0$-estimate of approximate solutions 
and the regularity of a weak solution, 
we will use rather the given orbifold metric $g$ itself 
than an {\it asymptotically} cylindrical metric $\bar{g} \in [g|_X]$ 
on $X := M - \{\check{p}\}$ 
with $\bar{g} = r^{-2}\cdot g$ near the singularity $\check{p}$,  
where $r(\cdot) := {\rm dist}_g(\cdot , \check{p})$.  

First, note that 
$$ 
Y(M, [g]_{orb}) = \inf_{u \in C^{\infty}_c(X), u \not\equiv 0} Q_{(X, g)}(u). 
$$ 
Let $B_{\rho}$ be the open geodesic ball centered at $\check{p}$ of radius $\rho > 0$ with respect to $g$. 
Set 
$$ 
Y_i := \inf_{u \in C^{\infty}_c(X-\overline{B_{1/i}}), u \not\equiv 0} Q_{(X, g)}(u) 
$$ 
for $i \in \mathbb{N}$. 
We have that 
$$ 
Y_i > Y_{i+1} > Y_{i+2} > \cdots, \qquad \qquad \qquad \qquad \ \ 
$$ 
$$  
\lim_{i\to \infty} Y_i = \inf_{u \in C^{\infty}_c(X), u \not\equiv 0} Q_{(X, g)}(u) 
= Y(M, [g]_{orb}). 
$$ 
It then follows from the strict inequality~(1) and the above 
that there exists a large integer $i_0$ such that 
$$ 
Y_i < Y(S^n)/|\Gamma|^{2/n} < Y(S^n) \quad {\rm for\ any}\ \ i \geq i_0. 
$$ 
Similarly to the case for compact manifolds without boundary, 
this implies that there exists a non-negative $Q_{(X - B_{1/i}, g)}$-minimizer 
$u_i \in C^{\infty}(X - B_{1/i})$ such that, for each $i \geq i_0$,  
$$ 
Q_{(X - B_{1/i}, g)}(u_i) = Y_i,\qquad \int_{X - B_{1/i}} u_i^{\frac{2n}{n-2}} d\mu_g = 1, 
$$ 
$$ 
u_i = 0\quad {\rm on}\ \ \partial B_{1/i},\qquad u_i > 0\quad {\rm in}\ \ X-\overline{B_{1/i}}. 
$$ 
We denote the zero extension of each $u_i$ to $M$ by also the same symbol $u_i$. 

Suppose that the sequence $\{u_i\}$ has a uniform $C^0$-bound, that is, 
there exists a constant $L > 0$ such that 
$$ 
||u_i||_{C^0(M)} \leq L\quad {\rm for}\ \ i \geq i_0. 
$$ 
Under this uniform $C^0$-estimate, 
then there exists a non-negative $Q_{(M, g)}$-minimizer $u \in W^{1,2}(M; g)$ with $||u||_{C^0(M)} \leq L$  
such that (taking a subsequence if necessary) 
$$ 
u_i \rightarrow u\quad {\rm weakly\ in}\ \ W^{1,2}(M; g),\qquad 
u_i \rightarrow u\quad {\rm strongly\ in}\ \ L^2(M; g).  
$$ 
Lebesgue's bounded convergence theorem combined with the above uniform $C^0$-estimate for $\{u_i\}$ implies that 
$$ 
\int_M u^{\frac{2n}{n-2}} d\mu_g = 1. 
$$ 
By this equation and the fact that $\{u_i\}$ is a $Q_{(M, g)}$-minimizing sequence, 
we have 
$$ 
u_i \rightarrow u\quad {\rm strongly\ in}\ \ W^{1,2}(M; g). 
$$ 
Under the $C^0$-estimate $||u||_{C^0(M)} \leq L$, 
applying the standard elliptic $L^p$-estimates to the Euler-Lagrange equations 
for $u$ on $X$ and the lifting $(\varphi^{-1}\circ \pi)^{\ast} u$ on $\mathbb{B}_{\tau}({\bf 0})$,  
we obtain that $u \in C^{\infty}(M)$. 
Here, $u \in C^{\infty}(M)$ means that $u \in C^{\infty}(X)$ 
and the lifting $(\varphi^{-1}\circ \pi)^{\ast} u$ is smooth on $\mathbb{B}_{\tau}({\bf 0})$. 
The maximum principle \cite[Proposition~3.75]{Aubin-Book} implies that $u > 0$ everywhere on $M$,  
and then we get an orbifold Yamabe metric  
$$ 
\tilde{g} := u^{4/(n-2)}\cdot g \in [g]_{orb}. 
$$ 

To complete the proof, we need only to show a uniform $C^0$-estimate for the sequence $\{u_i\}$. 
For each $u_i$, take a maximum point $q_i \in X$ of $u_i$, 
and set $m_i := u_i(q_i)$. 
Taking a subsequence if necessary, 
we then have that there exists a point $q_{\infty} \in M$ such that 
$$ 
\lim_{i \to \infty} q_i = q_{\infty}. 
$$ 
Suppose that 
$$ 
\lim_{i \to \infty} m_i = \infty.  
$$ 
Then, we will lead to a contradiction as below. \\ 
\underline{Case~1.}\ \ $q_{\infty} \ne \check{p}$~: 
Let $\{V, x = (x^1, \cdots, x^n)\}$ be a geodesic normal coordinate system 
centered at $q_{\infty}$ satisfying $V \subset X$. 
We may assume that $\{ |x| < 1 \} \subset V$. 
Set 
$$ 
v_i(x) := m_i^{-1}\cdot u_i\big{(}m_i^{-\frac{2}{n-2}}\cdot x + x(q_i)\big{)}\quad 
{\rm for}\ \ x \in \{ |x| < m_i^{\frac{2}{n-2}}(1 - |x(q_i)|) \}. 
$$ 
Similarly to the proof of Theorem~2.1 in \cite[Chapter~5]{ScYa-Book}, 
there exists a positive function $v \in C^{\infty}(\mathbb{R}^n)$ such that 
$$ 
v_i \rightarrow v\quad {\rm in\ the}~C^2
\textrm{-topology\ on\ each\ relatively\ compact\ domain\ in}\ \ \mathbb{R}^n.  
$$ 
Hence, $v$ satisfies the following: 
$$ 
- \alpha_n \Delta_0 v = Y(M, [g]_{orb})\cdot v^{\frac{n+2}{n-2}}\quad {\rm on}\ \ \mathbb{R}^n, 
$$ 
$$ 
\int_{\mathbb{R}^n} v^{\frac{2n}{n-2}} dx \leq 
\liminf_{i \to \infty} \int_V u_i^{\frac{2n}{n-2}} d\mu_g \leq 1, 
$$ 
where $\Delta_0$ denotes the Laplacian with respect to the Euclidean metric. 
This implies that $Y(M, [g]_{orb}) \geq Y(S^n)$, 
and then it contradicts to the assumption~(1). \\ 
\underline{Case~2.}\ \ $q_{\infty} = \check{p}$~: 
In this case, we consider rather the liftings $\tilde{u}_i := (\varphi^{-1}\circ \pi)^{\ast} u_i$ 
on $\mathbb{B}_{\tau}({\bf 0})$ than $u_i$ themselves. 
Similarly to the above, set 
$$ 
\tilde{v}_i(x) := m_i^{-1}\cdot \tilde{u}_i\big{(}m_i^{-\frac{2}{n-2}}\cdot x + x(q_i)\big{)}\quad 
{\rm for}\ \ x \in \{ x \in \mathbb{R}^n~|~|x| < m_i^{\frac{2}{n-2}}(\tau - |x(q_i)|) \}. 
$$ 
Then, there exists a positive function $\tilde{v} \in C^{\infty}(\mathbb{R}^n)$ such that 
$$ 
\tilde{v}_i \rightarrow \tilde{v}\quad {\rm in\ the}~C^2
\textrm{-topology\ on\ each\ relatively\ compact\ domain\ in}\ \ \mathbb{R}^n.  
$$ 
Moreover, $\tilde{v}$ satisfies the following: 
$$ 
- \alpha_n \Delta_0 \tilde{v} = Y(M, [g]_{orb})\cdot \tilde{v}^{\frac{n+2}{n-2}}\quad 
{\rm on}\ \ \mathbb{R}^n, 
$$ 
$$ 
\int_{\mathbb{R}^n} \tilde{v}^{\frac{2n}{n-2}} dx \leq 
\liminf_{i \to \infty} \int_{\mathbb{B}_{\tau}({\bf 0})} \tilde{u}_i^{\frac{2n}{n-2}} d\mu_{\hat{g}} \leq |\Gamma|,  
$$ 
where $\hat{g} := (\varphi^{-1}\circ \pi)^{\ast} g$. 
This implies that 
$$ 
Y(M, [g]_{orb}) \geq \frac{Y(S^n)}{|\Gamma|^{2/n}}, 
$$  
and then it also contradicts to the assumption~(1). 
\qquad \qquad \qquad \qquad \qquad \qquad \qquad $\square$ \\ 

We can now prove the second assertion in Theorem~A. \\ 
\quad \\ 
{\it Proof of the second assertion in Theorem~A}.\ \ 
First, we note that 
$$ 
Y(M_1, [\check{g}]_{orb}) \leq Y^{orb}(M_1) \leq Y(N) \leq 0. 
$$ 
It then follows from Theorem~3.1 and the above inequality 
that there exists a constant scalar curvature orbifold metric 
$\tilde{g} \in [\check{g}]_{orb}$ satisfying 
$$ 
E(\tilde{g}) = Y(M_1, [\check{g}]_{orb}) \leq 0. 
$$ 
Similarly to the case for smooth conformal manifolds, 
the uniqueness of constant scalar curvature orbifold metrics 
in a non-positive orbifold conformal class \cite[Lemma~2.3]{AkBo-2} 
implies that, up to a scaling,  
$$ 
\check{g} = \tilde{g}. 
$$ 
Combining the above with the assumption $E(\check{g}) = Y(N)$,   
we then have 
$$ 
Y(N) = E(\check{g}) = Y(M_1, [\check{g}]_{orb}) \leq Y^{orb}(M_1) \leq Y(N). 
$$ 
This implies that $Y^{orb}(M_1) = Y(M_1, [\check{g}]_{orb}) = Y(N)$.  
\qquad \qquad \qquad \qquad \qquad \qquad \ \ $\square$\\

\section{\bf Exact computations}\label{EC} 

We first prove Theorem~B.\\ 
\quad \\ 
{\it Proof of Theorem~B.}\ \ 
First, note that the canonical model $X_{can}$ is obtained by blowing down 
each connected component of the union of the $(-2)$-curves in $X$ into a point. 
The structure of an open neighborhood of each singular point in $X_{can}$ 
is modeled by one of A-D-E singularities, 
that is, the quotient singularity $\mathbb{C}^2/\Gamma$ 
with a non-trivial finite subgroup $\Gamma < SU(2)$. 
Then, $X_{can}$ admits a K\"ahler-Einstein orbifold metric $\check{g}$ \cite{Kor} 
satisfying 
$$ 
E(\check{g}) = - 4 \sqrt{2} \pi \sqrt{2 \chi(X) + 3 \tau(X)}.  
$$ 
We denote the singularities of $X_{can}$ by 
$\{(\check{p}_1, \Gamma_1), \cdots, (\check{p}_{\ell}, \Gamma_{\ell})\}$. 

For each $\Gamma_j$, let $X_j$ denote the minimal resolution of $\mathbb{C}^2/\Gamma_j$. 
Then, each $X_j$ admits a hyperK\"ahler ALE metric $h_j$ \cite{Kr}, 
and $(X_j, h_j)$ has a smooth conformal compactification 
$(\hat{X_j} := X_j \sqcup \{\infty_j\}, \hat{h}_j)$ 
with singularity $\{(\infty_j, \Gamma_j)\}$ \cite{ChLeWe, Vi}, 
which has a positive Yamabe constant $Y(\hat{X}_j, [\hat{h}_j]_{orb}) > 0$. 

With these understandings, 
$X$ can be decomposed by 
$$ 
X = X_{can} \#_{\sqcup_{j=1}^{\ell} (S^3/\Gamma_j)} 
\big{(}\sqcup_{j = 1}^{\ell} \hat{X}_j \big{)}.  
$$ 
By Theorem~A, this combined with $Y(X) < 0$ and $Y^{orb}(\hat{X}_j) > 0$ implies  
$$ 
Y^{orb}(X_{can}) \leq Y(X) < 0. 
$$ 
Recall that the K\"ahler-Einstein orbifold metric $\check{g}$ satisfies 
$$ 
E(\check{g}) = - 4 \sqrt{2} \pi \sqrt{2 \chi(X) + 3 \tau(X)} = Y(X).   
$$ 
This gives the desired conclusion:  
$$ 
\qquad \qquad \qquad \qquad \quad 
Y^{orb}(X_{can}) = Y(X_{can}, [\check{g}]_{orb}) = Y(X).
\qquad \qquad \qquad \qquad \quad \square\\ 
$$ 

Finally, we give two more typical exact computations of the orbifold Yamabe invariant. \\ 
\quad \\ 
\underline{\bf 1.}\quad Let $T$ be a complex $2$-dimensional torus 
and $\check{T} := T/\langle{\rm id}, \iota\rangle$ the quotient $4$-orbifold with 
$16$-singularities $\{(\check{p}_1, \langle{\rm id}, \iota\rangle), \cdots, (\check{p}_{16}, \langle{\rm id}, \iota\rangle)\}$. 
Here, $\langle{\rm id}, \iota\rangle~(\cong \mathbb{Z}_2)$ denotes the group of degree $2$ generated by 
$$ 
\iota : \mathbb{C}^2 \rightarrow \mathbb{C}^2,\ \ (z_1, z_2) \mapsto (- z_1, - z_2). 
$$ 
Pushing down the flat metric on $T$ to $\check{T}$, 
we have a flat orbifold metric $\check{g}_{flat}$ on $\check{T}$. 

\begin{prop} 
$$ 
Y^{orb}(\check{T}) = Y(\check{T}, [\check{g}_{flat}]_{orb}) = 0. 
$$ 
\end{prop} 
\noindent  
{\it Proof.}\ \ Let $\mathcal{O}(-2)$ denote the complex line bundle 
over the complex projective line $\mathbb{C}P^1$ of degree $- 2$. 
Then, there exists a cylindrical metric $\bar{g}$ on $\mathcal{O}(-2)$ 
modeled by $(S^3/\langle{\rm id}, \iota\rangle, h_0)$ with positive scalar curvature $R_{\bar{g}} > 0$ 
(cf.~\cite[Example~4.1.27]{Ni-Book}). 
Hence, $(\mathcal{O}(-2), \bar{g})$ has a smooth conformal 
compactification $(\widehat{\mathcal{O}(-2)} := \mathcal{O}(-2) \sqcup \{\infty\}, \hat{g})$ 
with singularity $\{(\infty, \langle{\rm id}, \iota\rangle)\}$. 
Note that, from the uniform positivity of $R_{\bar{g}}$ and the Sobolev embedding 
$W^{1, 2}(\mathcal{O}(- 2); \bar{g}) 
\hookrightarrow L^4(\mathcal{O}(- 2); \bar{g})$, 
$$  
Y(\widehat{\mathcal{O}(-2)}, [\hat{g}]_{orb}) 
= Y(\mathcal{O}(-2), [\bar{g}]) > 0. 
$$  
Let $(N_1, H_1), \cdots, (N_{16}, H_{16})$ be 
the $16$-copies of $(\widehat{\mathcal{O}(-2)}, \langle{\rm id}, \iota\rangle)$. 
With these understandings, 
the generalized connected sum 
$$ 
X := \check{T} \#_{\sqcup_{j=1}^{\ell} (S^3/H_j)} (\sqcup_{j=1}^{\ell} N_j) 
$$ 
is diffeomorphic to the Kummer surface, and hence $Y(X) = 0$. 
By Theorem~A, we then have 
$$ 
Y^{orb}(\check{T}) \leq Y(X) = 0. 
$$ 
Note that 
$$ 
E(\check{g}_{flat}) = 0 = Y(X), 
$$ 
and hence 
$$ 
\qquad \qquad \qquad \qquad \qquad 
Y^{orb}(\check{T}) = Y(\check{T}, [\check{g}_{flat}]_{orb}) 
= Y(X) = 0. \qquad \qquad \qquad \qquad \square 
$$ 
\quad \\ 
\underline{\bf 2.}\quad 
Let $\Sigma$ be an exotic sphere of dimension 
$n := 8k + 2 \geq 10$ with $\alpha([\Sigma]) \ne 0$, 
where $\alpha$ is the $\alpha$-homomorphism 
from the spin cobordism group $\Omega^{spin}_n$ 
to the $KO$-group $KO^{-n}(pt) \cong \mathbb{Z}_2$ (cf.~\cite[Chapter~2]{LaMi-Book}).  
For any integer $\ell \geq 2$, 
set 
$$ 
G_{\ell} := \{ \zeta^j I \in GL(4k+1; \mathbb{C})~|~j = 0, \cdots, \ell -1 \},\quad 
\zeta := \exp(2\pi \sqrt{-1}/\ell) \in \mathbb{C},  
$$ 
where $I$ denotes the identity matrix.  
The finite group $G_{\ell}$ acts the $n$-sphere 
$S^n \subset \mathbb{R}^{n+1} = \mathbb{C}^{4k+1} \times \mathbb{R}$ 
by   
$$ 
A : \mathbb{C}^{4k+1} \times \mathbb{R} \rightarrow \mathbb{C}^{4k+1} \times \mathbb{R},\ \ 
(z, t) \mapsto A\cdot (z, t) := (A\cdot z, t)\quad {\rm for}\ \ A \in G_{\ell}.
$$ 
Then, the quotient space $S^n/G_{\ell}$ is a compact $n$-orbifold 
with two singularities 
$\{(\check{p}_+ := [(0, \cdots, 0, 1)], G^+_{\ell} := G_{\ell}), 
(\check{p}_- := [(0, \cdots, 0, - 1)], G^-_{\ell} := G_{\ell})\}$. 
Pushing down the standard metric $g_0$ on $S^n$ to $S^n/G_{\ell}$, 
we have an orbifold metric $\check{g}_0$ of constant curvature one on $S^n/G_{\ell}$. 
Note that the space $((S^n/G_{\ell}) - \{\check{p}_+, \check{p}_-\}, \check{g}_0)$ 
is conformal to the product space 
$((S^{n-1}/G_{\ell}) \times \mathbb{R}, \bar{g} := h_0 + dt^2)$. 
Then, this combined with $R_{\bar{g}} = R_{h_0} = (n-1)(n-2) > 0$ and 
the Sobolev embedding 
$W^{1, 2}((S^{n-1}/G_{\ell}) \times~\mathbb{R}; \bar{g}) 
\hookrightarrow L^{2n/(n-2)}((S^{n-1}/G_{\ell}) \times \mathbb{R}; \bar{g})$ 
implies that 
\begin{equation} 
Y^{orb}(S^n/G_{\ell}) \geq Y(S^n/G_{\ell}, [\check{g}_0]_{orb}) 
= Y((S^{n-1}/G_{\ell}) \times \mathbb{R}, [\bar{g}]) > 0. 
\end{equation}  

\begin{prop} 
$$ 
Y^{orb}(\Sigma \# (S^n/G_{\ell})) = 0. 
$$ 
Here, $\Sigma \# (S^n/G_{\ell})$ stands for the connected sum 
of $\Sigma$ and $S^n/G_{\ell}$ in the usual sense. 
\end{prop} 
\noindent  
{\it Proof.}\ \ 
We first note the following. 
By results of Lichnerowicz and Hitchin (cf.~\cite[Chapters~2, 4]{LaMi-Book}) 
for $\alpha : \Omega^{spin}_n \rightarrow \mathbb{Z}_2$, 
$Y(\Sigma) \leq 0$. 
On the other hand, 
Petean \cite{Petean-2} proved that 
any simply connected compact manifold of dimension greater than $4$ 
has a non-negative Yamabe invariant. 
Hence, we have 
$$ 
Y(\Sigma) = 0. 
$$ 

Let 
$$ 
N_{\ell} := (S^n/G_{\ell}) \#_{(S^{n-1}/G^+_{\ell}) \sqcup (S^{n-1}/G^-_{\ell})} (\overline{S^n/G_{\ell}})  
$$ 
denotes the generalized connected sum.  
Here, $\overline{S^n/G_{\ell}}$ is the same $n$-orbifold, 
but equipped with the opposite orientation. 
It turns out that 
$$ 
N_{\ell} = (S^{n-1}/G_{\ell}) \times S^1, 
$$ 
and then it is a compact spin $n$-manifold with positive Yamabe invariant. 
Then, the positivity $Y(N_{\ell})~>~0$ implies that $\alpha([N_{\ell}]) = 0$, 
and hence 
$\alpha([\Sigma \# N_{\ell}]) = \alpha([\Sigma]) + \alpha([N_{\ell}]) \ne 0$. 
Therefore, 
\begin{equation} 
Y(\Sigma \# N_{\ell}) = 0. 
\end{equation} 
 
We now decompose $\Sigma \# N_{\ell}$ as the generalized connected sum 
$$ 
\Sigma \# N_{\ell} = \big{(}\Sigma \# (S^n/G_{\ell})\big{)} 
\#_{(S^{n-1}/G^+_{\ell}) \sqcup (S^{n-1}/G^-_{\ell})} (\overline{S^n/G_{\ell}}). 
$$ 
It then follows from Theorem~A combined with (2), (3) that 
\begin{equation} 
Y^{orb}(\Sigma \# (S^n/G_{\ell})) \leq Y(\Sigma \# N_{\ell}) = 0. 
\end{equation} 
On the other hand, Kobayashi's inequality for $Y^{orb}(\Sigma \# (S^n/G_{\ell}))$ 
still holds. 
Hence, 
\begin{equation}
0 = Y(\Sigma) = \min \{Y(\Sigma), Y^{orb}(S^n/G_{\ell}) \} 
\leq Y^{orb}(\Sigma \# (S^n/G_{\ell})).  
\end{equation} 
The inequalities (4), (5) give the desired conclusion: 
$$ 
\qquad \qquad \qquad \qquad \qquad \qquad \qquad 
Y^{orb}(\Sigma \# (S^n/G_{\ell})) = 0. 
\qquad \qquad \qquad \qquad \qquad \quad \square
$$ 
\quad \\ 
\quad \\

\bibliographystyle{amsbook}

\vspace{20mm} 

\end{document}